\documentclass[12pt]{article}

\usepackage{graphicx}
\usepackage[francais]{babel}
\usepackage{amsfonts}
\usepackage{pstricks}
\usepackage{amsmath}
\usepackage{amssymb}
\usepackage{amsthm}
\usepackage{empheq}
\usepackage{indentfirst}
\input{amssym.def}
\input{amssym}
\allowdisplaybreaks

\renewcommand{\theequation}{\arabic{section}.\arabic{equation}}

\newcommand{\tbd}{\tilde{\bd}}

\newcommand{\tbu}{\tilde{\bu}}
\newcommand{\tr}{\tilde{\rho}}
\newcommand{\bu}{\mathbf{u}}

\newcommand{\bd}{\mathbf{d}}
\newcommand{\bof}{\mathbf{f}}
\newcommand{\bs}{\mathbb{S}}

\newcommand{\R}{\mathbb{R}}

\newcommand{\pa}{\mathbb{\partial}}

\def\dis{\displaystyle}

\newtheorem{theorem}{\bf Theorem}[section]

\newtheorem{remark}{\bf Remark}[section]

\def\md{\mathrm{d}}

\title{Weak-strong uniqueness property\\ for the compressible flow of liquid crystals}

\author{Yong-Fu Yang,$^{1,2,}$\thanks{Corresponding author.}\,\, Changsheng Dou,$^{2}$ and  Qiangchang Ju$^2$ }

\date{}
\begin{document}
\maketitle

\vspace{-3mm}

\begin{center}
{\small $^1$ Department of Mathematics, College of Sciences, Hohai University,\\[2mm]
 Nanjing 210098, Jiangsu Province, P.R. China  \\[2mm]
$^2$ Institute of Applied Physics and Computational Mathematics,\\[2mm]
 P.O. Box 8009, Beijing 100088, P.R. China\\[2mm]
Email : fudanyoung@gmail.com, \hspace{2mm}douchangsheng@163.com,
\hspace{2mm}qiangchang\_ju@yahoo.com}
\end{center}

% --------------------------------------------------------------------------

\vspace{2mm}
%\textcolor[rgb]{1.00,0.00,0.00}{}
\begin{center}
\begin{minipage}{14cm}
        {\bf Abstract.} {\small Weak-strong uniqueness property in the class of finite energy weak solutions is established for two different compressible liquid crystal systems by the method of relative entropy. To overcome the difficulties caused by the molecular direction with inhomogeneous Dirichlet boundary condition, new techniques are introduced to build up the relative entropy inequalities.}
\end{minipage}
\end{center}

\vspace{3mm}


%$1%%%%%%%%%%%%%%%%%%%%%%%%%%%%%%%%%%%%%%%%%%%%%%%%%%%%%%%%%%%%%%%%%%%%%%

\section{Introduction}
\newcounter{intro}
\renewcommand{\theequation}{\thesection.\theintro}

In this paper, we consider the equations for the compressible flow of liquid crystals. We are first concerned with the following simplified version of Ericksen-Leslie system (for physical background see \cite{Eri62,Les68,Lin89,Lin95,LL95,LL96})~:
\stepcounter{intro}
\begin{equation}
\label{lcd-1}
\left\{\begin{array}{lll}
 \dis   \rho_t + \text{div}(\rho\, \mathbf{u}) = 0, \\[2mm]
 \dis   (\rho\, \mathbf{u})_t + \text{div}(\rho \mathbf{u}\otimes \mathbf{u})+\nabla P = \mu \triangle \bu - \lambda \,\text{div}\Big(\nabla \bd \odot \nabla \bd -\big(\frac{1}{2}|\nabla \bd|^2+F(\nabla \bd)\mathbb{I}_3\big)\Big),\\[3mm]
 \dis   \bd_t+\bu \cdot \nabla \bd= \theta (\triangle \bd-\textbf{f}(\bd)),
\end{array}\right.
\end{equation}
where
$\rho \geq 0$ is the density, $\mathbf{u}=(u_1,u_2,u_3)$ is the velocity, $\mathbf{d}=(d_1,d_2,d_3)$ is the molecular direction field. $\mu ,\lambda,\theta$ are positive constants. $P(\rho)$ is the pressure-density
function and here we consider the case
\begin{equation}
P(\rho)=a \rho^\gamma,
\label{aad}
\end{equation}
where  $a>0 $ and $ \gamma >1$ are constants.
 $\mathbb{I}_3$ is the $3 \times 3$ unit matrix. The term $\nabla \bd \odot \nabla \bd$ is defined by
\[
\nabla \bd \odot \nabla \bd = (\nabla \bd)^T\nabla \bd,
\]
where $(\nabla \bd)^T$ stands for the transpose of the $3 \times 3$ matrix $\nabla \bd$. The vector-valued smooth function $\bof(\bd)$ and the scalar function $F(\bd)$ are assumed here to satisfy the relation
\[
\dis \bof(\bd)=\nabla_{\bd}F(\bd).
\]
For example, we can choose $ F(\bd)$ and  $\bof(\bd)$ as
\[
\dis F(\bd) = \frac{1}{4\sigma_0^2}(|\bd|^2-1)^2,\quad \bof(\bd) = \frac{1}{\sigma_0^2}(|\bd|^2-1)\bd,
\]
where $F(\bd)$ is the Ginzburg-Landau penalization and $\sigma_0 > 0$ a constant.

Let $\Omega \in \R^3$ be a bounded smooth domain. In this paper, we will consider the following initial-boundary conditions~:
\stepcounter{intro}
\begin{equation}\label{ini-1}
    \dis (\rho,\rho \bu,\bd)|_{t=0} = (\rho_0(x), \mathbf{m}_0(x),\bd_0(x)), \quad x \in \Omega,
\end{equation}
and
\stepcounter{intro}
\begin{equation}\label{bou-1}
    \dis \bu|_{\pa \Omega}=0, \quad \bd|_{\pa \Omega} = \bd_0(x) \quad x \in \pa\Omega,
\end{equation}
where
\[
\dis \rho_0 \in L^{\gamma}(\Omega), \quad \rho_0 \geq 0; \quad \bd_0 \in L^{\infty}(\Omega)\cap H^1(\Omega);
\]
\[
\dis \mathbf{m}_0 \in L^1(\Omega),\quad \mathbf{m}_0=0 \,\, \text{if}\,\,\rho_0 = 0; \quad \frac{|\mathbf{m}_0|^2}{\rho_0} \in L^1(\Omega).
\]

The global existence of weak solutions $(\rho,\bu,\bd)$ to the initial-boundary problem (\ref{lcd-1})-(\ref{bou-1})in three-dimension with $\dis \gamma >3/2$ was obtained by Wang and Yu in \cite{WY12} and Liu and Qing in \cite{LQ12}, independently. In addition, when system (\ref{lcd-1}) is incompressible, Jiang and Tan \cite{JT09} and Liu and Zhang in \cite{LZ09} proved the global weak existence of solutions to the flow of nematic liquid crystals for fluids with non-constant density. Based on the existence result, Dai {\it et al.} in \cite{DQS12} extended the regularity and uniqueness results of Lin and Liu in \cite{LL95} to the systems of nematic liquid crystals with non-constant fluid density. In this paper, we are interested in the uniqueness of the weak solution obtained in \cite{WY12} and \cite{LQ12}. To our best knowledge, so far there are very few results concerning uniqueness of weak solutions to the initial-boundary value problem (\ref{lcd-1})--(\ref{bou-1}).

Second, we consider the following the hydrodynamic flow equation of nematic compressible liquid crystals:
\stepcounter{intro}
\begin{equation}
\label{lcd-2}
\left\{\begin{array}{lll}
 \dis   \rho_t + \text{div}(\rho\, \mathbf{u}) = 0, \\[2mm]
 \dis   (\rho\, \mathbf{u})_t + \text{div}(\rho \mathbf{u}\otimes \mathbf{u})+\nabla P = \mu \triangle \bu - \lambda \,\text{div}(\nabla \bd \odot \nabla \bd -\frac{1}{2}|\nabla \bd|^2\mathbb{I}_3),\\[3mm]
 \dis   \bd_t+\bu \cdot \nabla \bd= \theta (\triangle \bd+|\nabla \bd|^2\bd),
\end{array}\right.
\end{equation}
where $\bd \in \bs^2$ and the other symbols have the same meanings with those in system (\ref{lcd-1}). In this situation, the Ossen-Frank energy configuration functional reduces to the Dirichlet energy functional. We refer to the readers to consult the recent papers \cite{DHWZ11} and \cite{LQ12} for the derivation of the system (\ref{lcd-2}). For the system (\ref{lcd-2}), we are concerned with the same initial conditions (\ref{ini-1}) but with the following different boundary conditions~:
\stepcounter{intro}
\begin{equation}\label{bou-2}
    \dis \bu|_{\pa \Omega}=0, \quad \frac{\pa \bd}{\pa \nu}|_{\pa \Omega} = 0 \quad x \in \pa\Omega,
\end{equation}
where $\nu$ is the unit outer normal vector of $\pa \Omega$.

In contrast with system (\ref{lcd-1}), from the mathematical point of view, it is much more difficult to deal with the nonlinear term $|\nabla \bd|^2 \bd$ appearing in the third equation of (\ref{lcd-2}). Even for the incompressible flow, there have been no satisfactory results concerning the global existence of weak solutions. Recently, Lin {\it et al.} in \cite{LLW10} proved the existence of global-in-time weak solutions on a bounded smooth domain in $\R^2$.
For three dimensional case, the problem is still open. We should mention that Li and Wang in \cite{LW11} has established its weak-strong uniqueness principle in three dimension, provided that the existence of its weak solution is obtained.

The compressible flow (\ref{lcd-2}) of liquid crystals is much more complicated and hard to study mathematically due to the compressibility. For the one-dimensional case, the global existence of smooth and weak solutions to the compressible flow of liquid crystals was obtained by Ding {\it et al.} in \cite{DWW11,DLWW12}. As for three dimensional case, Huang {\it et al.} in \cite{HWW12} and Liu and Zhang in \cite{LZ09} established the local existence of a unique strong solution provided that the initial data are sufficiently regular and satisfy a natural compatibility condition. However, the global existence of weak solution to the compressible flow of liquid crystals in multi-dimension is still open. In this paper, we shall establish the weak-strong uniqueness property for (\ref{lcd-2}).

Recently,  Feireisl, Jin and Novotn\'y \cite{FJN12} established the weak-strong uniqueness for the compressible Navier-Stokes equations. They showed that a weak solution coincides with the strong solution with the same initial data, if the strong solution exists. We also refer to the readers to consult the recent papers \cite{FNS11,FN12,Ge06,Ge11} for more weak-strong uniqueness results. Motivated by Feireisl {\it
et al.} \cite{FJN12}, we aim to establish the weak-strong uniqueness property for two simplified Erichsen-Leslie system (\ref{lcd-1}) and (\ref{lcd-2}), respectively. Our method is essentially based on the relative entropy,
the modified relative entropy inequality, and a Gronwall-type argument.
Compared with the existence result of weak solution in \cite{LQ12,WY12}, where both of them require $\dis \gamma > 3/2$, we are going to make use of the techniques, established in \cite{YDJ12}, to estimate the remainder $\mathcal{R}$ (or $\mathcal{R}_1$ (for the definition see (\ref{R}) (or (\ref{R1}))), so as to establish the weak-strong uniqueness property for an improved lower bound for any adiabatic exponent $\gamma >1$. When dealing with the compressible nematic liquid crystal flow (\ref{lcd-1}) or (\ref{lcd-2}), the main difficulty lies in the coupling and interaction between the velocity filed $\bu$ and the direction filed $\bd$. In particular, we should emphasize here that, as far as the weak-strong uniqueness property for the initial-boundary problem (\ref{lcd-1})-(\ref{bou-1}) is concerned, more efforts are required to build up the relative entropy inequality so as to overcome the effects arising from inhomogeneous boundary condition $\bd|_{\pa \Omega}=\bd_0(x)$, which is quite different from the case of the isentropic compressible Navier-stokes investigated by Feireisl {\it et al.} in \cite{FJN12}. According to the definition of weak solutions for the compressible flow of liquid crystals, it is easy to see that
$$\dis \bd \in L^2(0,T;H^2(\Omega)),$$
so the weak solution $\bd$ actually solves the third equation of (\ref{lcd-1}) (or (\ref{lcd-2})) in the strong sense.
Consequently, we can combine the arguments of Feireisl {\it et al.} in \cite{FJN12} with the energy estimates for parabolic equation to obtain the desired relative entropy inequality.

The sizes of the positive constants $\mu$, $\lambda$, and $\theta$ do not play important roles in our proofs, we shall therefore assume, for simplicity, that
\stepcounter{intro}
\begin{equation}\label{norm}
    \dis \mu = \lambda = \theta =1,
\end{equation}
throughout this paper. In addition, to simplify the notations, we always set
$$\textbf{A} \cdot (\nabla \textbf{B})\textbf{C} = ((\textbf{C} \cdot \nabla)\textbf{B})\cdot \textbf{A},$$
in which $\textbf{A}$, $\textbf{B}$, and $\textbf{C}$ are vectors in $\R^3$.

The paper of the rest is organized as follows. In the next section, we recall
the definitions of weak and strong solutions to the compressible flow of liquid crystals for two kinds of models and state the main results. Section 3 is devoted to the derivation of the relative entropy inequality and the proof of Theorem \ref{T2.1}. Finally, we prove Theorem \ref{T2.2} in Section 4.

%$2%%%%%%%%%%%%%%%%%%%%%%%%%%%%%%%%%%%%%%%%%%%%%%%%%%%%%%%%%%%%%%%%%%%%%%

\section{Main results}
\newcounter{result}
\renewcommand{\theequation}{\thesection.\theresult}

$\{\rho,\bu,\bd\}$ is a
{\it finite energy weak solution} (see \cite{Lions98,NS04}) to the initial-boundary value problem
(\ref{lcd-1})--(\ref{bou-1}), if, for any $T>0$,

$\bullet$ $\rho \geq 0,$ \quad $\rho \in L^{\infty}(0,T;L^{\gamma}(\Omega))$,\quad $\bu \in L^2(0,T; H^1_0(\Omega)),$
\[
\dis \bd \in L^{\infty}((0,T)\times \Omega)\cap L^{\infty}(0,T;H^1(\Omega)) \cap L^2(0,T;H^2(\Omega)),
\]
\par \quad with $(\rho,\rho\bu,\bd)(0,x)=(\rho_0(x),\mathbf{m}_0(x),\bd_0(x))$ for $x \in \Omega$;

$\bullet$ The first equation in (\ref{lcd-1}) is replaced by a family of integral identities
\stepcounter{result}
\begin{equation}
\label{cont}
\dis \int_{\Omega}\,\rho(\tau,\cdot)\varphi(\tau,\cdot)\,\md x - \int_{\Omega}\,\rho_0\varphi(0,\cdot)\,\md x = \int_0^{\tau}\int_{\Omega}\,
(\rho\partial_t\varphi+\rho\mathbf{u}\cdot \nabla \varphi) \,\md x \md t
\end{equation}
\par $\quad$for any $\varphi \in C^1([0,T]\times \bar{\Omega})$, and any $\tau \in [0,T]$;

$\bullet$ Momentum equations $(\ref{lcd-1})_2$ are satisfied in the sense of distributions, specifically,
\stepcounter{result}
\begin{align}
\label{momen}
\dis & \int_{\Omega}\,\rho\mathbf{u}(\tau,\cdot)\varphi(\tau,\cdot)\,\md x - \int_{\Omega}\,\rho_0\mathbf{u}_0\cdot\varphi(0,\cdot)\,\md x \\\nonumber
\dis =& \int_0^{\tau}\int_{\Omega}\,
\Big(\rho\mathbf{u}\cdot\partial_t\varphi+\rho\mathbf{u}\otimes\mathbf{u}: \nabla \varphi + P(\rho) \text{div} \varphi-\nabla \mathbf{u}:\nabla \varphi\Big)\,\md x\md t\\\nonumber
&\,+\int_0^{\tau}\int_{\Omega}\,\Big(\nabla \bd \odot \nabla \bd - \big(\frac{1}{2}|\nabla \bd|^2+F(\bd)\big)\mathbb{I}_3\Big):\nabla \varphi \,\md x \md t
\end{align}
\par $\quad$for any $\varphi \in C^1([0,T]\times\bar{\Omega})$, $\varphi|_{\pa \Omega}=0$, and any $\tau \in [0,T]$.

$\bullet$ Equations (\ref{lcd-1})$_3$ are
replaced by a family of integral identities
\stepcounter{result}
\begin{align}
\label{mag} \dis\quad \quad\quad &\int_{\Omega}\,
\bd(\tau,\cdot)\cdot \varphi(\tau,\cdot)\,\md x -\int_{\Omega}\,
\bd_0\cdot \varphi(0,\cdot)\,\md x\\\nonumber =&
\int_0^{\tau}\int_{\Omega}\,\Big( \bd \cdot \pa_t
\varphi-\nabla\bd:\nabla \varphi - \varphi\cdot(\nabla
\bd)\bu-\varphi\cdot \bof(\bd)\Big)\, \md x \md t
\end{align}
\par $\quad$for any $\varphi \in C^1([0,T]\times\bar{\Omega})$, $\varphi|_{\pa \Omega}=0$, and any $\tau \in [0,T]$.

$\bullet$ The energy inequality
\stepcounter{result}
\begin{equation}
\label{enegy-1}
\dis E(t)+\int_0^{\tau}\int_{\Omega}\,\Big( |\nabla \bu|^2+|\triangle \bd - \bof(\bd)|^2
\Big)\,\md x \md t \leq E(0)
\end{equation}
\par $\quad$ holds for a.e. $\tau \in [0,T]$, where
$$ \dis E(t)=\int_{\Omega}\,\Big(\frac{1}{2}\rho|\bu|^2+\frac{a}{\gamma-1}\rho^{\gamma}+\frac{1}{2}|\nabla \bd|^2+F(\bd)\Big)\,\md x,$$
\par $\quad$ and
$$ \dis E(0)=\int_{\Omega}\,\Big(\frac{1}{2}\frac{|\mathbf{m}_0|^2}{\rho_0}+\frac{a}{\gamma-1}\rho_0^{\gamma}+\frac{1}{2}|\nabla \bd_0|^2+F(\bd_0)\Big)\,\md x.$$

The existence of global-in-time finite energy weak solutions to the initial-boundary problem (\ref{lcd-1})-(\ref{bou-1}) with the adiabatic exponent $\gamma > \frac{3}{2}$ was established in \cite{LQ12,WY12}, provided there exists a positive constant $C_0$ such that $\bd \cdot \bof(\bd) \geq 0$ for all $|\bd| \geq C_0 >0$.

$\{\tr,\tbu,\tbd\}$ is called a classical (strong) solution to the initial-boundary problem (\ref{lcd-1})-(\ref{bou-1}) in $(0,T)\times \Omega$ if
\stepcounter{result}
\begin{equation}
\label{regu}
\left\{
\begin{array}{ll}
\dis \tr\in C^1([0,T]\times \bar{\Omega}),\quad \tr(t,x) \geq \underline{\rho}>0
\quad \text{for all } (t,x)\in (0,T)\times \Omega ,\\[2mm]
\dis \tbu,\,\pa_t \tbu,\,\nabla^2\tbu\in
C([0,T]\times\bar{\Omega}),\,\tbd,\,\pa_t \tbd,\, \triangle\tbd\in
C([0,T]\times\bar{\Omega})
\end{array}
\right.
\end{equation}
and $\tr, \tbu, \tbd$ satisfy equation (\ref{lcd-1}), together with
the boundary conditions (\ref{bou-1}). Observe that hypothesis
(\ref{regu}) requires the following regularity properties of the
initial data:
\stepcounter{result}
\begin{equation}
\label{initial}
\left\{
\begin{array}{ll}
\dis \tr(0,\cdot)=\rho_{0}\in C^1(\bar{\Omega}),\quad \rho_0 \geq \underline{\rho}>0,\\[2mm]
\dis \tbu(0,\cdot)=\bu_0 \in C^2(\bar{\Omega}),\quad \tbd(0,\cdot)=\bd_0 \in C^2(\bar{\Omega}).
\end{array}
\right.
\end{equation}

We are now ready to state the first result of this paper.

\begin{theorem}
\label{T2.1}
Let $\Omega \in \R^3$ be a bounded domain with a boundary of class $C^{2+\kappa}$, $\kappa>0$, and $\gamma>1$. Suppose that $\{\rho,\bu,\bd\}$ is a finite energy weak solution to the initial-boundary problem (\ref{lcd-1})-(\ref{bou-1}) in $(0,T)\times\Omega$ in the sense specified above, and suppose that $\{\tr,\tbu,\tbd\}$ is a strong solution emanating from the same initial data (\ref{initial}).

Then
\[\dis \rho \equiv \tr,\quad \bu \equiv \tbu,\quad \bd \equiv \tbd.
\]
\end{theorem}

In a similar way, we define the {\it finite energy weak solution} $\{\rho_1,\bu_1,\bd_1\}$ to the initial-boundary value problem (\ref{lcd-2}), (\ref{ini-1}), and (\ref{bou-2}) in the following sense~: for any $T>0$,

$\bullet$ $\rho_1 \geq 0,$ \quad $\rho_1 \in L^{\infty}(0,T;L^{\gamma}(\Omega))$,\quad $\bu_1 \in L^2(0,T; H^1_0(\Omega)),$
\[
\dis \bd_1 \in L^{\infty}((0,T)\times \Omega)\cap L^{\infty}(0,T;H^1(\Omega)) \cap L^2(0,T;H^2(\Omega)),
\]
\par \quad with $(\rho_1,\rho_1\bu_1,\bd_1)(0,x)=(\rho_0(x),\mathbf{m}_0(x),\bd_0(x))$ for $x \in \Omega$;

$\bullet$ Similar to (\ref{cont})-(\ref{mag}), the equations (\ref{lcd-2}) hold in $\mathcal{D}'((0,T)\times \Omega)$

$\bullet$ The energy inequality
\stepcounter{result}
\begin{equation}
\label{enegy-2}
\dis E_1(t)+\int_0^{\tau}\int_{\Omega}\,\Big( |\nabla \bu_1|^2+\big|\triangle \bd_1 + |\nabla \bd_1|^2\bd_1\big|^2
\Big)\,\md x \md t \leq E_1(0)
\end{equation}
\par $\quad$ holds for a.e. $\tau \in [0,T]$, where
$$ \dis E_1(t)=\int_{\Omega}\,\Big(\frac{1}{2}\rho_1 |\bu_1|^2+\frac{a}{\gamma-1}\rho_1^{\gamma}+\frac{1}{2}|\nabla \bd_1|^2\Big)\,\md x,$$
\par $\quad$ and
$$ \dis E_1(0)=\int_{\Omega}\,\Big(\frac{1}{2}\frac{|\mathbf{m}_0|^2}{\rho_0}+\frac{a}{\gamma-1}\rho_0^{\gamma}+\frac{1}{2}|\nabla \bd_0|^2\Big)\,\md x.$$

\begin{remark}
\label{R2.2}
In the derivation of energy inequality (\ref{enegy-2}), we have used the fact that $|\bd_1|=1$ to get
\[
\dis \Big(\pa_t \bd_1 + \bu_1 \cdot \nabla \bd_1\Big) \cdot |\nabla \bd_1|^2 \bd_1 = \frac{1}{2}|\nabla \bd_1|^2 \Big(\pa_t |\bd_1|^2 + \bu_1 \cdot \nabla |\bd_1|^2\Big) = 0.
\]
\end{remark}

\begin{remark}
\label{R2.3}
We should remark that the global-in time renormalized finite energy weak solutions to the initial-boundary value problem (\ref{lcd-2}), (\ref{ini-1}) and (\ref{bou-2}) is still open. For one-dimensional case, Ding {\it et al.} in \cite{DWW11,DLWW12} obtained the global-in-time existence of weak solutions.
\end{remark}

Similarly, $\{\tr_1,\tbu_1,\tbd_1\}$ is called a classical (strong) solution to the initial-boundary value problem (\ref{lcd-2}), (\ref{ini-1}) and (\ref{bou-2}) in $(0,T)\times \Omega$ if
\stepcounter{result}
\begin{equation}
\label{regu-2}
\left\{
\begin{array}{ll}
\dis \tr_1\in C^1([0,T]\times \bar{\Omega}),\quad \tr_1(t,x) \geq \underline{\rho}>0
\quad \text{for all } (t,x)\in (0,T)\times \Omega ,\\[2mm]
\dis \tbu_1,\,\pa_t \tbu_1,\,\nabla^2\tbu_1\in
C([0,T]\times\bar{\Omega}),\,\tbd_1,\,\pa_t \tbd_1,\, \triangle\tbd_1\in
C([0,T]\times\bar{\Omega})
\end{array}
\right.
\end{equation}
and $\tr_1, \tbu_1, \tbd_1$ satisfy equation (\ref{lcd-2}), together with
the boundary conditions (\ref{bou-2}). Observe that hypothesis
(\ref{regu-2}) requires the following regularity properties of the
initial data:
\stepcounter{result}
\begin{equation}
\label{initial-2}
\left\{
\begin{array}{ll}
\dis \tr_1(0,\cdot)=\rho_{0}\in C^1(\bar{\Omega}),\quad \rho_0 \geq \underline{\rho}>0,\\[2mm]
\dis \tbu_1(0,\cdot)=\bu_0 \in C^2(\bar{\Omega}),\quad \tbd_1(0,\cdot)=\bd_0 \in C^2(\bar{\Omega}).
\end{array}
\right.
\end{equation}

We now end up this section with another result of this paper.

\begin{theorem}
\label{T2.2}
Let $\Omega \in \R^3$ be a bounded domain with a boundary of class $C^{2+\kappa}$, $\kappa>0$, and $\gamma>1$. Suppose that $\{\rho_1,\bu_1,\bd_1\}$ is a finite energy weak solution to the initial-boundary value problem (\ref{lcd-2}), (\ref{ini-1}) and (\ref{bou-2}) in $(0,T)\times\Omega$ in the sense specified above, and suppose that $\{\tr_1,\tbu_1,\tbd_1\}$ is a strong solution emanating from the same initial data (\ref{initial-2}).

Then
\[\dis \rho_1 \equiv \tr_1,\quad \bu_1 \equiv \tbu_1,\quad \bd_1 \equiv \tbd_1.
\]
\end{theorem}

%$3%%%%%%%%%%%%%%%%%%%%%%%%%%%%%%%%%%%%%%%%%%%%%%%%%%%%%%%%%%%%%%%%%%%%%%
\section{Proof of Theorem \ref{T2.1}}
\newcounter{ineq}
\renewcommand{\theequation}{\thesection.\theineq}

Motivated by the concept of {\it relative entropy} in \cite{FJN12}, we first define {\it relative entropy}
%%%%
$\dis \mathcal{E} = \mathcal{E}\left([\rho,\bu,\bd]|[\tilde{\rho},\tilde{\bu},\tilde{\bd}]\right)$, with respect to $\{\tilde{\rho},\tilde{\bu},\tilde{\bd}\}$, as
\stepcounter{ineq}
\begin{equation}
\label{re}
\dis \mathcal{E} = \int_{\Omega} \,\big( \frac{1}{2}\rho|\mathbf{u}-\tilde{\mathbf{u}}|^2+\Pi(\rho)-\Pi'(\tilde{\rho})(\rho-\tilde{\rho})-\Pi(\tilde{\rho})
+\frac{1}{2}|\nabla \bd-\nabla\tilde{ \bd}|^2\big)\, \md x,
\end{equation}
where
\stepcounter{ineq}
\begin{equation}
\label{pr}
\dis \Pi(\rho)=\frac{a}{\gamma-1} \rho^{\gamma}.
%\rho\int_{\bar{\rho}}^{\rho}\,\frac{p(z)}{z^2}\, \md z,\quad \quad \bar{\rho}=\frac{1}{|\Omega|}\int_{\Omega}\,\rho\, \md x.
\end{equation}

In this section, we are going to deduce a relative entropy inequality satisfied by any weak solution to the initial-boundary value problem (\ref{lcd-1})-(\ref{bou-1}). To this end, consider a triple $\{\tr,\tbu,\tbd\}$ of smooth functions, $\tr$ bounded away from zero in $[0,T]\times \bar{\Omega}$, $\tbu|_{\pa \Omega}=0$, and $\tbd|_{\pa \Omega}=\bd_0(x)$. In addition, $\tbu$ and $\tbd$ solve the third equation of (\ref{lcd-1}).

Noticing that the boundary conditions for $\bd$ and $\tbd$ are inhomogeneous, i.e. $\bd|_{\pa \Omega}=\tbd|_{\pa \Omega}=\bd_0(x)$, we should adapt and modify the arguments in \cite{FJN12} to build up the relative entropy inequality. More regularity of $\bd$ allows us to make use of energy estimates for parabolic equation, which help us overcome the difficulty due to the inhomogeneous boundary conditions. Consequently, combining the arguments in \cite{FJN12} with the energy estimates for parabolic equations yields the desired relative entropy inequality.

To begin with, we take $\tbu$ as a test function in the momentum equation (\ref{momen}) to obtain
\stepcounter{ineq}
\begin{align}
\label{momen-1}
\dis &\,\int_{\Omega}\,\rho\bu\cdot \tbu(\tau,\cdot)\,\md x -\int_{\Omega}\, \rho_0\bu_0\cdot \tbu(0,\cdot)\,\md x\\\nonumber
= &\,\int_0^{\tau}\int_{\Omega}\Big(\rho\bu\cdot \pa_t \tbu + \rho \bu \otimes \bu : \nabla \tbu + P(\rho) \text{div}\tbu-\nabla \bu: \nabla \tbu \Big)\, \md x \md t\\\nonumber
&\,- \int_0^{\tau}\int_{\Omega}\Big((\triangle \bd- \bof(\bd)) \cdot (\nabla \bd) \tbu\Big)\, \md x \md t.
\end{align}

Second, we can use the scalar quantity $\dis \varphi=\frac{1}{2}|\tbu|^2$ and $\varphi=\Pi'(\tr)$, respectively, as test functions in the continuity equation (\ref{cont}) to get
\stepcounter{ineq}
\begin{equation}
\label{cont-1}
\dis \int_{\Omega}\,\frac{1}{2}\rho|\tbu|^2(\tau,\cdot)\,\md x= \int_{\Omega}\,\frac{1}{2}\rho_0|\tbu|^2(0,\cdot)\,\md x + \int_0^{\tau}\int_{\Omega}\,\Big(\rho \tbu \cdot \pa_t \tbu+\rho \tbu \cdot (\nabla \tbu)\bu\Big)\,\md x \md t
\end{equation}
and
\stepcounter{ineq}
\begin{equation}
\label{cont-2}
\dis \int_{\Omega}\,\rho\Pi'(\tr)(\tau,\cdot)\, \md x=\int_{\Omega}\,\rho_0 \Pi(\tr)(0,\cdot)\,\md x + \int_0^{\tau}\int_{\Omega}\,\Big(\rho \pa_t \Pi'(\tr)+\rho \bu \cdot \nabla \Pi'(\tr) \Big)\,\md x \md t.
\end{equation}
Since $(\bu,\bd)$ solves the third equation of (\ref{lcd-1}) in the strong sense, it is easy to see that
\[
\dis \pa_t (\bd-\tbd)+\bu \cdot \nabla \bd - \tbu \cdot \nabla \tbd = (\triangle \bd - \triangle \tbd)- (\bof(\bd) -\bof(\tbd)), \quad \quad a.e..
\]
Multiply the above equation by $(\triangle \bd - \triangle \tbd)$ and integrate over $\Omega \times (0,\tau)$, we have
\stepcounter{ineq}
\begin{align}\label{d-1}
    \dis \,&\,\int_{\Omega}\,\frac{1}{2}|\nabla \bd - \nabla \tbd|^2 (\tau,\cdot)\, \md x + \int_0^{\tau}\int_{\Omega}\,|\triangle \bd - \triangle \tbd|^2\, \md x \md t\\\nonumber
    = \,&\, \int_{\Omega}\,\frac{1}{2}|\nabla \bd_0 - \nabla \tbd(0,\cdot)|^2 \, \md x + \int_0^{\tau}\int_{\Omega}\, (\triangle \bd - \triangle \tbd)\cdot (\nabla \bd)\bu \, \md x \md t \\\nonumber
    \,&\,\,-\int_0^{\tau}\int_{\Omega}\, (\triangle \bd - \triangle \tbd)\cdot (\nabla \tbd)\tbu \, \md x \md t + \int_0^{\tau}\int_{\Omega}\, (\bof(\bd) - \bof(\tbd))\cdot (\triangle \bd - \triangle \tbd) \md x \md t.
\end{align}
Finally, multiplying the third equation of (\ref{lcd-1}) by $\triangle \bd - \bof(\bd)$ gives
\stepcounter{ineq}
\begin{align}\label{d-2}
    \dis \,&\,\int_{\Omega}\,(\frac{1}{2} |\nabla \bd|^2+F(\bd))(\tau,\cdot)\,\md x + \int_0^{\tau}\int_{\Omega}\,|\triangle \bd - \bof(\bd)|^2\, \md x \md t \\\nonumber
    =\,&\, \int_{\Omega}\,\big(\frac{1}{2} |\nabla \bd_0|^2+F(\bd_0)\big)\,\md x + \int_0^{\tau}\int_{\Omega}\,(\triangle \bd - \bof(\bd)) \cdot (\nabla \bd) \bu \, \md x \md t.
\end{align}

Summing up the relations (\ref{momen-1})-(\ref{d-1}) with the energy inequality (\ref{enegy-1}), we infer that
\stepcounter{ineq}
\begin{align}
\label{entr}
\dis \, &\,\int_{\Omega}\,\Big(\frac{1}{2}\rho|\bu-\tbu|^2+\Pi(\rho)-\rho\Pi'(\tr)+\frac{1}{2}|\nabla \bd-\nabla \tbd|^2\Big)(\tau,\cdot)\,\md x\\\nonumber
&\,+\int_0^{\tau}\int_{\Omega}\,|\nabla \bu-\nabla \tbu|^2\, \md x \md t + \int_0^{\tau}\int_{\Omega}\, |\triangle \bd - \triangle \tbd|^2\, \md x \md t\\\nonumber
\leq \,&\, \int_{\Omega}\,\Big(\frac{1}{2}\rho_0|\bu_0-\tbu(0,\cdot)|^2+\Pi(\rho_0)-\rho_0\Pi'(\tr(0,\cdot))+\frac{1}{2}|\nabla \bd_0-\nabla \tbd(0,\cdot)|^2\Big)\,\md x\\\nonumber
& \,+ \int_0^{\tau}\int_{\Omega}\,\rho(\pa_t \tbu + \bu \cdot \nabla \tbu)\cdot(\tbu-\bu)\,\md x \md t + \int_0^{\tau}\int_{\Omega}\, \nabla \tbu : \nabla (\tbu-\bu) \, \md x \md t\\\nonumber
&\, -\int_0^{\tau}\int_{\Omega}\,(\rho \pa_t \Pi'(\tr)+\rho \bu \cdot \nabla \Pi'(\tr))\,\md x \md t - \int_0^{\tau}\int_{\Omega}\, P(\rho) \text{div} \tbu \, \md x \md t\\\nonumber
&\,- \int_0^{\tau}\int_{\Omega}\, (\triangle \bd-\bof(\bd)) \cdot (\nabla \bd)\bu\, \md x \md t + \int_0^{\tau}\int_{\Omega}\, (\triangle \bd-\bof(\bd)) \cdot (\nabla \bd)\tbu\, \md x \md t\\\nonumber
&\,+\int_0^{\tau}\int_{\Omega}\, (\triangle \bd-\triangle \tbd) \cdot \big((\nabla \bd)\bu - (\nabla \tbd)\tbu\big)\, \md x \md t \\\nonumber
&\,+\int_0^{\tau}\int_{\Omega}\, (\triangle \bd-\triangle \tbd) \cdot \big(\bof(\bd)-\bof(\tbd)\big)\, \md x \md t,
\end{align}
where we have used (\ref{d-2}). By virtue of the definition (\ref{pr}) of $\Pi$, it is easy to see that
\begin{equation*}
    \dis \Pi'(\tr)\tr -\Pi(\tr) = P(\tr)
\end{equation*}
and
\begin{equation*}
    \dis \int_{\Omega}\,\Big(\tr \pa_t \Pi'(\tr)+\tr \nabla \Pi'(\tr)\cdot \tbu + P(\tr) \text{div}\tbu\Big)\, \md x = \int_{\Omega}\,\pa_t P(\tr)\,\md x.
\end{equation*}
As a consequence, we deduce from the identities
\begin{equation*}
    \dis \int_{\Omega}\,P(\tr)(\tau,\cdot)\,\md x - \int_{\Omega}\,P(\tr)(0,\cdot)\, \md x
    = \int_0^{\tau}\int_{\Omega}\,\pa_t P(\tr)\,\md x \md t
\end{equation*}
and (\ref{entr}) that the desired relative entropy inequality holds~:
\stepcounter{ineq}
\begin{align}
\label{entr-1}
    \dis \,&\, \mathcal{E}(\tau) + \int_0^{\tau}\int_{\Omega}\,|\nabla \bu- \nabla \tbu|^2\, \md x \md t + \int_0^{\tau}\int_{\Omega}\, |\triangle \bd - \triangle \tbd|^2\, \md x \md t\\\nonumber
\leq \,&\, \mathcal{E}(0) + \int_0^{\tau}\, \mathcal{R}(\rho,\bu,\bd,\tr,\tbu,\tbd)\,\md t,
\end{align}
where
\stepcounter{ineq}
\begin{equation}\label{R}
    \dis \mathcal{R} = \mathcal{R}(\rho,\bu,\bd,\tr,\tbu,\tbd):= \mathcal{R}_d + \mathcal{R}_c,
\end{equation}
\stepcounter{ineq}
\begin{align}\label{Rd}
    \dis \mathcal{R}_d:=\,&\,\int_{\Omega}\,\rho (\tbu-\bu)\cdot \big(\pa_t \tbu + (\nabla \tbu)\bu\big) \, \md x +\int_{\Omega}\, \nabla \tbu: \nabla(\tbu-\bu) \, \md x \\\nonumber
    &\,+\int_{\Omega}\,\big((\tr-\rho)\pa_t \Pi'(\tr) + \nabla \Pi'(\tr) \cdot (\tr \tbu- \rho \bu)\big)\, \md x - \int_{\Omega}\,\text{div}\tbu \big(P(\rho)-P(\tr)\big)\,\md x,
\end{align}
and
\stepcounter{ineq}
\begin{align}\label{Rc}
    \dis \mathcal{R}_c:=\,&\,- \int_0^{\tau}\int_{\Omega}\, (\triangle \bd-\bof(\bd)) \cdot (\nabla \bd)\bu\, \md x \md t + \int_0^{\tau}\int_{\Omega}\, (\triangle \bd-\bof(\bd)) \cdot (\nabla \bd)\tbu\, \md x \md t\\\nonumber
&\,+\int_0^{\tau}\int_{\Omega}\, (\triangle \bd-\triangle \tbd) \cdot \big((\nabla \bd)\bu - (\nabla \tbd)\tbu\big)\, \md x \md t \\ \nonumber &\,+\int_0^{\tau}\int_{\Omega}\, (\triangle \bd-\triangle \tbd) \cdot \big(\bof(\bd)-\bof(\tbd)\big)\, \md x \md t.
\end{align}

In what follows, we shall finish the proof of Theorem \ref{T2.1} by applying the relative entropy inequality
(\ref{entr-1}) to $\{\tr,\tbu,\tbd\}$, where $\{\tr,\tbu,\tbd\}$ is a classical (smooth) solution of the initial-boundary value problem (\ref{lcd-1}), (\ref{initial}), and (\ref{bou-1}), such that
\begin{equation*}
    \dis \tr(0,\cdot)=\rho_0, \quad \tbu(0,\cdot)=\bu_0, \quad \tbd(0,\cdot)=\bd_0.
\end{equation*}

Accordingly, the integrals depending on the initial data on the right-hand side of (\ref{entr-1}) vanish,
and we apply a Gronwall type argument to deduce the desired result, namely,
\begin{equation*}
    \dis \rho \equiv \tr,\quad \bu \equiv \tbu,\quad \bd \equiv \tbd.
\end{equation*}
Our purpose is to examine all terms in the remainder (\ref{R}) and to show that they can be "absorbed"
by the left-hand side of (\ref{entr-1}).

Compared with \cite{FJN12}, we should remark that the main difficulty comes from the coupling and interaction between the velocity field $\bu$ and the direction field $\bd$. Moreover, in the context of the weak-strong uniqueness,
we only assume the adiabatic exponent $\gamma >1$, but not $\gamma >3/2$ as in \cite{FJN12}.

Similar to the proof of Theorem 2.1 in \cite{YDJ12} (see also \cite{FJN12}), we use (\ref{Rd}) to find that
\stepcounter{ineq}
\begin{align}\label{Rd-1}
    \dis \mathcal{R}_d = \,&\, \int_{\Omega}\, \rho (\tbu-\bu)\cdot \big((\nabla \tbu) (\bu-\tbu)\big) \, \md x + \int_{\Omega}\, \frac{(\rho-\tr)}{\tr} \triangle \tbu \cdot (\tbu-\bu)\,\md x\\\nonumber
    \,&\,-\int_{\Omega}\,\text{div}\tbu\Big(P(\rho)-P'(\tr)(\rho-\tr)-P(\tr)\Big)\, \md x \\\nonumber
    \,&\,- \int_{\Omega}\,\frac{\rho}{\tr}\text{div}\Big(\nabla \tbd \odot \nabla \tbd - \big(\frac{1}{2}|\nabla \tbd|^2+F(\tbd)\big) \mathbb{I}_3  \Big)\cdot (\tbu-\bu) \, \md x \\\nonumber
    = \,&\, \int_{\Omega}\, \rho (\tbu-\bu)\cdot \big((\nabla \tbu) (\bu-\tbu)\big) \, \md x \\\nonumber \,&\,-\int_{\Omega}\,\text{div}\tbu\Big(P(\rho)-P'(\tr)(\rho-\tr)-P(\tr)\Big)\, \md x \\\nonumber
    \,&\,+ \int_{\Omega}\, \frac{(\rho-\tr)}{\tr}\left[\triangle \tbu- \text{div}\Big(\nabla \tbd \odot \nabla \tbd - \big(\frac{1}{2}|\nabla \tbd|^2+F(\tbd)\big) \mathbb{I}_3\Big)\right]\cdot (\tbu-\bu)\,\md x\\\nonumber
    \,&\,- \int_{\Omega}\,(\triangle \tbd- \bof(\tbd)) \cdot (\nabla \tbd) (\tbu-\bu) \, \md x\\\nonumber
    =: &\, \overline{\mathcal{R}}_d - \int_{\Omega}\,(\triangle \tbd- \bof(\tbd)) \cdot (\nabla \tbd) (\tbu-\bu) \, \md x.
\end{align}
After a tedious but straightforward computation, it follows from (\ref{Rc}) that
\stepcounter{ineq}
\begin{align}\label{Rc-1}
   \overline{\mathcal{R}}_c := \,&\, \mathcal{R}_c - \int_{\Omega}\,(\triangle \tbd- \bof(\tbd)) \cdot (\nabla \tbd) (\tbu-\bu) \, \md x\\\nonumber
    = \,\, & \, \int_{\Omega}\,  (\triangle \bd-\triangle \tbd)\cdot (\bof (\bd) - \bof (\tbd))\,\md x + \int_{\Omega}\,(\triangle \bd-\triangle \tbd)\cdot (\nabla \bd - \nabla \tbd) \tbu \, \md x \\\nonumber
    \,&\, + \int_{\Omega}\, \triangle \tbd \cdot (\nabla \bd - \nabla \tbd)(\tbu - \bu)\, \md x + \int_{\Omega}\, \bof(\bd) \cdot (\nabla \bd - \nabla \tbd) (\bu - \tbu)\, \md x\\\nonumber
    \,&\, + \int_{\Omega}\,(\bof(\bd)-\bof(\tbd)) \cdot \nabla \tbd (\bu - \tbu)\, \md x.
\end{align}
Consequently, it follows from the definitions of $\overline{\mathcal{R}}_d$ and $\overline{\mathcal{R}}_c$ that
\stepcounter{ineq}
\begin{equation}\label{R'}
    \dis \mathcal{R}(\rho,\bu,\bd,\tr,\tbu,\tbd) = \overline{\mathcal{R}}_d + \overline{\mathcal{R}}_c.
\end{equation}

In order to prove the weak-strong uniqueness property, we have to estimate the remainder $\mathcal{R}$. As for $\overline{\mathcal{R}}_d$.
Since the procedures are almost same as that in \cite{YDJ12} (see also \cite{FJN12}), we just follow \cite{YDJ12}, list the outlines, and skip the details. In the sequel, we are going to focus on the estimation of $\overline{\mathcal{R}}_c$.

From (\ref{initial}), it is clear to see that
\stepcounter{ineq}
\begin{align}\label{a1}
    \dis \, & \, \left|\int_{\Omega}\, \rho (\tbu-\bu)\cdot \big((\nabla \tbu) (\bu-\tbu)\big) \, \md x - \int_{\Omega}\,\text{div}\tbu\Big(P(\rho)-P'(\tr)(\rho-\tr)-P(\tr)\Big)\, \md x\right| \\\nonumber
    \leq \,&\, C \|\nabla \tbu \|_{L^{\infty}(\Omega)}\, \mathcal{E}\Big([\rho,\bu,\bd]|[\tr,\tbu,\tbd]\Big).
\end{align}
Here and hereafter $C$ stands for a generic constant, which may change from line to line.

Let
$$\dis \tilde{\mathbf{g}}=\tilde{\mathbf{g}}(\triangle \tbu,\tbd,\nabla \tbd,\triangle \tbd)=\triangle \tbu- \text{div}\Big(\nabla \tbd \odot \nabla \tbd - \big(\frac{1}{2}|\nabla \tbd|^2+F(\tbd)\big) \mathbb{I}_3\Big).$$
Obviously, we have
\stepcounter{ineq}
\begin{align}\label{a3}
   \dis \,&\, \int_{\Omega}\, \frac{1}{\tr}(\rho-\tr)\, \tilde{\mathbf{g}} \cdot (\tbu-\bu)\,\md x \\\nonumber
    = \,&\,\int_{\{\frac{\tr}{2} < \rho < 2\tr\}}\, \frac{1}{\tr}(\rho-\tr)\,\tilde{\mathbf{g}} \cdot (\tbu-\bu)\,\md x + \int_{\{0 \leq \rho < \frac{\tr}{2}\}}\, \frac{1}{\tr}(\rho-\tr)\, \tilde{\mathbf{g}} \cdot (\tbu-\bu)\,\md x \\\nonumber
     \,&\,+\int_{\{ \rho \geq 2\tr\}}\, \frac{1}{\tr}(\rho-\tr)\, \tilde{\mathbf{g}} \cdot (\tbu-\bu)\,\md x.
\end{align}
Similar to that in \cite{YDJ12} (see also \cite{FJN12}), we make use of H\"{o}lder's inequality and Sobolev's inequality to show that
\stepcounter{ineq}
\begin{align}\label{a4}
    \dis \,&\,\left| \left(\int_{\{\frac{\tr}{2} < \rho < 2\tr\}}+\int_{\{0 \leq \rho < \frac{\tr}{2}\}}\right)\, \frac{1}{\tr}(\rho-\tr)\, \tilde{\mathbf{g}} \cdot (\tbu-\bu)\,\md x\right| \\\nonumber
    \leq \,&\, C(\delta)\|\frac{\tilde{\mathbf{g}}}{\tr}\|^2_{L^3(\Omega)}\,\mathcal{E}\Big([\rho,\bu,\bd]|[\tr,\tbu,\tbd]\Big) + \delta \, \|\nabla \bu -\nabla \tbu\|_{L^2(\Omega)}^2
\end{align}
for any $\delta > 0$. On the other hand, noticing that
\[
\dis \Pi(\rho) - \Pi'(\tr)(\rho-\tr) - \Pi(\tr) \geq C \rho^{\gamma},\quad \quad \text{as}\quad \rho \geq 2 \tr \geq 2 \underline{\rho}
\]
and
\[
\dis \left|\frac{\rho-\tr}{\rho \tr}\right|\,\rho^{\frac{1}{2}-\frac{\gamma}{2}} \leq C, \quad \quad \text{as}\quad \rho \geq 2 \tr \geq 2 \underline{\rho}\quad \text{and}\quad \gamma >1,
\]
we conclude from the definition of relative entropy $\mathcal{E}\Big([\rho,\bu,\bd]|[\tr,\tbu,\tbd]\Big)$ that
\stepcounter{ineq}
\begin{align}\label{a5}
    \dis \, & \, \left|\int_{\{ \rho \geq 2\tr\}}\, \frac{1}{\tr}(\rho-\tr)\,\tilde{\mathbf{g}} \cdot (\tbu-\bu)\,\md x \right|\\\nonumber
 \dis   = \, & \, \int_{\{ \rho \geq 2\tr\}}\, \left|\frac{\rho-\tr}{\rho \tr}\right|\,\rho^{\frac{1}{2}} \, \tilde{\mathbf{g}} \cdot \rho^{\frac{1}{2}}|\tbu-\bu|\,\md x\\\nonumber
\dis    = \, & \, \int_{\{ \rho \geq 2\tr\}}\,\left( \left|\frac{\rho-\tr}{\rho \tr}\right|\,\rho^{\frac{1}{2}-\frac{\gamma}{2}}\right)\rho^{\frac{\gamma}{2}} \, \tilde{\mathbf{g}} \cdot \rho^{\frac{1}{2}}|\tbu-\bu|\,\md x \\\nonumber
\dis    \leq \, & \, C \| \tilde{\mathbf{g}} \|_{L^{\infty}(\Omega)} \Big(\int_{\Omega}\, \rho^{\gamma}\, \md x\Big)^{\frac{1}{2}}\Big(\int_{\Omega}\, \frac{\rho}{2} |\bu-\tbu|^2\, \md x \Big)^{\frac{1}{2}}\\\nonumber
\dis    \leq \, & \, C \| \tilde{\mathbf{g}} \|_{L^{\infty}(\Omega)}\,\mathcal{E}\Big([\rho,\bu,\bd]|[\tr,\tbu,\tbd]\Big).
\end{align}

Next, we continue to estimate $\overline{\mathcal{R}}_c$.
Note that
\stepcounter{ineq}
\begin{equation}\label{a10}
    \|\bd\|_{L^{\infty}((0,T)\times \Omega)},\|\tbd\|_{L^{\infty}((0,T)\times \Omega)} \leq C
\end{equation}
and the fact that $\bof$ is smooth, we have
\stepcounter{ineq}
\begin{align}\label{a6}
    \dis \,&\, \left|\int_{\Omega}\,  (\triangle \bd-\triangle \tbd)\cdot (\bof (\bd) - \bof (\tbd))\,\md x\right| \\\nonumber
     \leq \,&\,C\,\left(\int_{\Omega}\,|\bd-\tbd|^2\,\md x\right)^{\frac{1}{2}} \left(\int_{\Omega}\,|\triangle \bd-\triangle \tbd|^2\,\md x\right)^{\frac{1}{2}}\\\nonumber
    \leq \,&\, \delta\, \|\triangle \bd - \triangle \tbd\|^2_{L^2(\Omega)}+ C(\delta)\, \|\bd - \tbd\|^2_{L^2(\Omega)}\\\nonumber
    \leq \,&\, \delta\, \|\triangle \bd - \triangle \tbd\|^2_{L^2(\Omega)}+ C(\delta)\, \|\nabla \bd - \nabla  \tbd\|^2_{L^2(\Omega)}\\\nonumber
    \leq \,&\, \delta\, \|\triangle \bd - \triangle \tbd\|^2_{L^2(\Omega)}+ C(\delta)\, \,\mathcal{E}\Big([\rho,\bu,\bd]|[\tr,\tbu,\tbd]\Big)
\end{align}
for any $\delta >0$, where we have used Sobolev's inequality. It follows from H\"{o}lder's inequality that
\stepcounter{ineq}
\begin{align}\label{a7}
    \dis \,&\, \left|\int_{\Omega}\,(\triangle \bd-\triangle \tbd)\cdot (\nabla \bd - \nabla \tbd) \tbu \, \md x  \right| \\\nonumber
   \leq  \,&\, \|\tbu\|_{L^{\infty}(\Omega)}\,\|\triangle \bd-\triangle \tbd\|_{L^2(\Omega)}\,\|\nabla \bd-\nabla \tbd\|_{L^2(\Omega)}\\\nonumber
   \leq  \,&\, \delta \, \|\triangle \bd-\triangle \tbd\|_{L^2(\Omega)}^2 + C(\delta)\, \|\tbu\|_{L^{\infty}(\Omega)}^2\|\nabla \bd-\nabla \tbd\|_{L^2(\Omega)}^2 \\\nonumber
   \leq  \,&\, \delta \, \|\triangle \bd-\triangle \tbd\|_{L^2(\Omega)}^2 + C(\delta)\, \|\tbu\|_{L^{\infty}(\Omega)}^2\,\mathcal{E}\Big([\rho,\bu,\bd]|[\tr,\tbu,\tbd]\Big)
\end{align}
for any $\delta >0$.
It is clear that
\stepcounter{ineq}
\begin{align}\label{a8}
    \dis \, & \, \left|\int_{\Omega}\, \triangle \tbd \cdot (\nabla \bd - \nabla \tbd)(\tbu - \bu)\, \md x + \int_{\Omega}\, \bof(\bd) \cdot (\nabla \bd - \nabla \tbd) (\bu - \tbu)\, \md x\right| \\\nonumber
    \leq \, & \, \big(\|\triangle \tbd\|_{L^{\infty}(\Omega)}+\|\bof (\bd)\|_{L^{\infty}(\Omega)}\big)\,\,\|\bu-\tbu\|_{L^2(\Omega)}\,\|\nabla \bd-\nabla \tbd\|_{L^2(\Omega)}\\\nonumber
    \leq \, & \, \delta\,\|\nabla \bu-\nabla \tbu\|_{L^2(\Omega)}^2 + C(\delta) \, \big(\|\triangle \tbd\|_{L^{\infty}(\Omega)}+\|\bof (\bd)\|_{L^{\infty}(\Omega)}\big)^2\,\,\|\nabla \bd-\nabla \tbd\|_{L^2(\Omega)}^2\\\nonumber
    \leq \, & \, \delta\,\|\nabla \bu-\nabla \tbu\|_{L^2(\Omega)}^2 + C(\delta) \, \big(\|\triangle \tbd\|_{L^{\infty}(\Omega)}+\|\bof (\bd)\|_{L^{\infty}(\Omega)}\big)^2\,\,\mathcal{E}\Big([\rho,\bu,\bd]|[\tr,\tbu,\tbd]\Big).
\end{align}
Finally, (\ref{a10}) and the fact that $\bof$ is smooth imply that
\stepcounter{ineq}
\begin{align}\label{a9}
    \dis \, & \, \left|\int_{\Omega}\, \big(\bof(\bd)-\bof(\tbd)\big) \cdot (\nabla \tbd)(\bu - \tbu)\, \md x\right| \\\nonumber
    \leq \, & \, C\, \|\nabla \tbd\|_{L^{\infty}(\Omega)}\,\|\bu-\tbu\|_{L^2(\Omega)}\,\|\bd-\tbd\|_{L^2(\Omega)}\\\nonumber
    \leq \, & \, \delta\,\|\nabla \bu-\nabla \tbu\|_{L^2(\Omega)}^2 + C(\delta) \, \big(\|\nabla \tbd\|_{L^{\infty}(\Omega)}\big)^2\,\,\|\nabla \bd-\nabla \tbd\|_{L^2(\Omega)}^2\\\nonumber
    \leq \, & \, \delta\,\|\nabla \bu-\nabla \tbu\|_{L^2(\Omega)}^2 + C(\delta) \, \big(\|\nabla \tbd\|_{L^{\infty}(\Omega)}\big)^2\,\,\mathcal{E}\Big([\rho,\bu,\bd]|[\tr,\tbu,\tbd]\Big).
\end{align}

Summing up relations (\ref{entr-1})--(\ref{a9}), we conclude that the relative entropy inequality yields the desired conclusion
\[
\dis \mathcal{E}\Big([\rho,\bu,\bd]|[\tr,\tbu,\tbd]\Big)(\tau) \leq \int_0^{\tau} \, h(t)
\mathcal{E}\Big([\rho,\bu,\bd]|[\tr,\tbu,\tbd]\Big)(t)\, \md t \quad \text{with some }\, h \in L^1(0,T).
\]
Thus, Theorem \ref{T2.1} immediately follows from Gronwall's inequality.

%$4%%%%%%%%%%%%%%%%%%%%%%%%%%%%%%%%%%%%%%%%%%%%%%%%%%%%%%%%%%%%%%%%%%%%%%
\section{Proof of Theorem \ref{T2.2}}
\newcounter{proof}
\renewcommand{\theequation}{\thesection.\theproof}

This section is devoted to proving the weak-strong uniqueness property for the initial-boundary value problem (\ref{lcd-2}), (\ref{ini-1}), and (\ref{bou-2}). Compared with the problem discussed in Section 3, we would like to point out two points~: (i) The term $|\nabla \bd_1|^2 \bd_1$ in the third equation of (\ref{lcd-2}) has higher nonlinearity than $\bof(\bd)$ in the third equation of (\ref{lcd-1}), so that more regularity results should be obtained to overcome the difficulty caused by $|\nabla \bd_1|^2 \bd_1$. (ii) Different from $\bd|_{\pa \Omega}=\bd_0(x)$, we now have homogeneous Neumann boundary condition for $\bd_1$ at hand, namely, $\dis \frac{\pa \bd_1}{\pa \nu}|_{\pa \Omega}=0$, which implies that a modified relative entropy is required
to prove the weak-strong uniqueness property for the initial-boundary value problem (\ref{lcd-2}), (\ref{ini-1}), and (\ref{bou-2}). To be precise, we define {\it relative entropy}
%%%%
$\dis \mathcal{E}_1 = \mathcal{E}_1\left([\rho_1,\bu_1,\bd_1]|[\tilde{\rho}_1,\tilde{\bu}_1,\tilde{\bd}_1]\right)$, with respect to $\{\tilde{\rho}_1,\tilde{\bu}_1,\tilde{\bd}_1\}$, as
\stepcounter{proof}
\begin{align}
\label{re-2}
\dis \mathcal{E}_1 = \int_{\Omega} \,\Big(\, &\, \frac{1}{2}\rho_1|\mathbf{u}_1-\tilde{\mathbf{u}}_1|^2+\big(\Pi(\rho_1)-\Pi'(\tilde{\rho}_1)(\rho_1-\tilde{\rho}_1)
-\Pi(\tilde{\rho_1})\big) \\\nonumber
\,&\,+ \frac{1}{2}\big(| \bd_1- \tilde{ \bd}_1|^2 + |\nabla \bd_1-\nabla\tilde{ \bd}_1|^2\big)\Big)\, \md x,
\end{align}
where $\Pi$ is defined by (\ref{pr}), $\{\rho_1,\bu_1,\bd_1\}$ is the finite energy weak solution to the initial-boundary value problem (\ref{lcd-2}), (\ref{ini-1}), and (\ref{bou-2}) in the sense of section 2, and $\{\tr_1,\tbu_1,\tbd_1\}$ are smooth functions, $\tr_1$ bounded away from zero in $[0,T]\times \bar{\Omega}$, $\tbu_1|_{\pa \Omega}=0$, and $\dis \frac{\pa \bd_1}{\pa \nu}|_{\pa \Omega}=0$. In addition, $\tbu_1$ and $\tbd_1$ solve the third equation of (\ref{lcd-2}).

We first establish the relative entropy inequality. For this purpose, take $\tbu_1$ as a test function in the momentum equation to obtain
\stepcounter{proof}
\begin{align}
\label{momen-2}
\dis &\,\int_{\Omega}\,\rho_1\bu_1\cdot \tbu_1(\tau,\cdot)\,\md x -\int_{\Omega}\, \rho_0\bu_0\cdot \tbu_1(0,\cdot)\,\md x\\\nonumber
= &\,\int_0^{\tau}\int_{\Omega}\Big(\rho_1\bu_1\cdot \pa_t \tbu_1 + \rho_1 \bu_1 \otimes \bu_1 : \nabla \tbu_1 + P(\rho_1) \text{div}\tbu_1-\nabla \bu_1: \nabla \tbu_1 \Big)\, \md x \md t\\\nonumber
&\,- \int_0^{\tau}\int_{\Omega}\triangle \bd_1 \cdot (\nabla \bd_1) \tbu_1\, \md x \md t.
\end{align}
Second, we can use the scalar quantity $\dis \frac{1}{2}|\tbu_1|^2$ and $\Pi'(\tr_1)$, respectively, as test functions in the continuity equation to get
\stepcounter{proof}
\begin{equation}
\label{cont-3}
\dis \int_{\Omega}\,\frac{1}{2}\rho_1|\tbu_1|^2(\tau,\cdot)\,\md x= \int_{\Omega}\,\frac{1}{2}\rho_0|\tbu_1|^2(0,\cdot)\,\md x + \int_0^{\tau}\int_{\Omega}\,\Big(\rho_1 \tbu_1 \cdot \pa_t \tbu_1+\rho_1 \tbu_1 \cdot (\nabla \tbu_1)\bu_1\Big)\,\md x \md t
\end{equation}
and
\stepcounter{proof}
\begin{equation}
\label{cont-4}
\dis \int_{\Omega}\,\rho_1\Pi'(\tr_1)(\tau,\cdot)\, \md x=\int_{\Omega}\,\rho_0 \Pi(\tr_1)(0,\cdot)\,\md x + \int_0^{\tau}\int_{\Omega}\,\Big(\rho_1 \pa_t \Pi'(\tr_1)+\rho_1 \bu_1 \cdot \nabla \Pi'(\tr_1) \Big)\,\md x \md t.
\end{equation}
 Notice that $(\bu_1,\bd_1)$ solves the third equation of (\ref{lcd-2}) in the strong sense, we have for almost everywhere
 %\stepcounter{proof}
 \begin{equation*}
    \dis \pa_t (\bd_1-\tbd_1)+\bu_1 \cdot \nabla \bd_1 - \tbu_1 \cdot \nabla \tbd_1 = (\triangle \bd_1 - \triangle \tbd_1) + (|\nabla \bd_1|^2 \bd_1 - |\nabla \tbd_1|^2 \tbd_1).
 \end{equation*}
 Multiplying the above equation by $(\bd_1-\tbd_1)$ and $(\triangle \bd_1 - \triangle \tbd_1)$, respectively, and integrating by parts, we obtain
\stepcounter{proof}
\begin{align}\label{d1}
    \dis \,&\, \int_{\Omega}\, \frac{1}{2} |\bd_1-\tbd_1|^2 \, \md x + \int_0^{\tau}\int_{\Omega}\, |\nabla \bd_1 - \nabla \tbd_1|^2 \, \md x \md t \\\nonumber
    = \,&\, \int_{\Omega}\, \frac{1}{2} |\bd_0-\tbd_1(0,\cdot)|^2 \, \md x + \int_0^{\tau}\int_{\Omega}\, \big(|\nabla \bd_1|^2 \bd_1 - |\nabla \tbd_1|^2 \tbd_1\big) \cdot (\bd_1 - \tbd_1) \, \md x \md t \\\nonumber
    \,&\, - \int_0^{\tau}\int_{\Omega}\, \big((\nabla \bd_1) \bu_1 - (\nabla \tbd_1) \tbu_1\big) \cdot (\bd_1 - \tbd_1) \, \md x \md t
\end{align}
and
\stepcounter{proof}
\begin{align}\label{d2}
    \dis \,&\, \int_{\Omega}\, \frac{1}{2} |\nabla \bd_1-\nabla \tbd_1|^2 \, \md x + \int_0^{\tau}\int_{\Omega}\, |\triangle \bd_1 - \triangle \tbd_1|^2 \, \md x \md t \\\nonumber
    = \,&\, \int_{\Omega}\, \frac{1}{2} |\nabla \bd_0-\nabla \tbd_1(0,\cdot)|^2 \, \md x - \int_0^{\tau}\int_{\Omega}\, \big(|\nabla \bd_1|^2 \bd_1 - |\nabla \tbd_1|^2 \tbd_1\big) \cdot (\triangle \bd_1 - \triangle \tbd_1) \, \md x \md t \\\nonumber
    \,&\, + \int_0^{\tau}\int_{\Omega}\, \big((\nabla \bd_1) \bu_1 - (\nabla \tbd_1) \tbu_1\big) \cdot (\triangle \bd_1 - \triangle \tbd_1) \, \md x \md t
\end{align}
Similarly, we multiply the third equation of (\ref{lcd-2}) by $\triangle \bd_1+|\nabla \bd_1|^2\bd_1$ to obtain
\stepcounter{proof}
\begin{align}\label{d3}
    \dis \,&\,\int_{\Omega}\, \frac{1}{2} |\nabla \bd_1|^2(\tau,\cdot)\,\md x + \int_0^{\tau}\int_{\Omega}\,\big|\triangle \bd_1 + |\nabla \bd_1|^2\bd_1\big|^2\, \md x \md t \\\nonumber
    =\,&\, \int_{\Omega}\,\frac{1}{2} |\nabla \bd_0|^2\,\md x + \int_0^{\tau}\int_{\Omega}\,\triangle \bd_1 \cdot (\nabla \bd_1) \bu_1 \, \md x \md t.
\end{align}

Summing up the relations (\ref{momen-2})-(\ref{d2}) with the energy inequality (\ref{enegy-2}), we infer that
\stepcounter{proof}
\begin{align}
\label{entr-2}
\dis \, &\,\int_{\Omega}\,\Big(\frac{1}{2}\rho|\bu_1-\tbu_1|^2+\Pi(\rho_1)-\rho_1\Pi'(\tr_1)+\frac{1}{2}\big(|\bd_1-\tbd_1|^2+|\nabla \bd_1-\nabla \tbd_1|^2\big)\Big)(\tau,\cdot)\,\md x\\\nonumber
&\,+\int_0^{\tau}\int_{\Omega}\,|\nabla \bu_1-\nabla \tbu_1|^2\, \md x \md t + \int_0^{\tau}\int_{\Omega}\, |\nabla \bd_1 - \nabla \tbd_1|^2\, \md x \md t + \int_0^{\tau}\int_{\Omega}\, |\triangle \bd_1 - \triangle \tbd_1|^2\, \md x \md t\\\nonumber
\leq \,&\, \int_{\Omega}\,\Big(\frac{1}{2}\rho_0|\bu_0-\tbu_1(0,\cdot)|^2+\Pi(\rho_0)-\rho_0\Pi'(\tr_1(0,\cdot))\\\nonumber
&\,+\frac{1}{2}\big(|\bd_0-\tbd_1(0,\cdot)|^2+|\nabla \bd_0-\nabla \tbd_1(0,\cdot)|^2\big)\Big)\,\md x\\\nonumber
& \,+ \int_0^{\tau}\int_{\Omega}\,\rho_1(\pa_t \tbu_1 + \bu_1 \cdot \nabla \tbu_1)\cdot(\tbu_1-\bu_1)\,\md x \md t + \int_0^{\tau}\int_{\Omega}\, \nabla \tbu_1 : \nabla (\tbu_1-\bu_1) \, \md x \md t\\\nonumber
&\, -\int_0^{\tau}\int_{\Omega}\,(\rho_1 \pa_t \Pi'(\tr_1)+\rho_1 \bu_1 \cdot \nabla \Pi'(\tr_1))\,\md x \md t - \int_0^{\tau}\int_{\Omega}\, P(\rho_1) \text{div} \tbu_1 \, \md x \md t\\\nonumber
&\,+ \int_0^{\tau}\int_{\Omega}\, \triangle \bd_1 \cdot (\nabla \bd_1)(\tbu_1-\bu_1)\, \md x \md t \\\nonumber
&\,+\int_0^{\tau}\int_{\Omega}\, (\triangle \bd_1-\triangle \tbd_1) \cdot \big((\nabla \bd_1)\bu_1 - (\nabla \tbd_1)\tbu_1\big)\, \md x \md t \\\nonumber
&\,-\int_0^{\tau}\int_{\Omega}\, (\triangle \bd_1-\triangle \tbd_1) \cdot \big(|\nabla \bd_1|^2\bd_1-|\nabla \tbd_1|^2\tbd_1\big)\, \md x \md t\\\nonumber
&\,+\int_0^{\tau}\int_{\Omega}\,(\bd_1 - \tbd_1) \cdot \big(|\nabla \bd_1|^2\bd_1-|\nabla \tbd_1|^2\tbd_1\big)\, \md x \md t\\\nonumber
&\,-\int_0^{\tau}\int_{\Omega}\, (\bd_1 - \tbd_1) \cdot \big((\nabla \bd_1)\bu_1 - (\nabla \tbd_1)\tbu_1\big)\, \md x \md t,
\end{align}
where we have used (\ref{d3}). Thus, we proceed the similar procedures as in the previous Section to obtain the desired relative entropy inequality
\stepcounter{proof}
\begin{align}
\nonumber
    \dis \,&\, \mathcal{E}_1(\tau) + \int_0^{\tau}\int_{\Omega}\,|\nabla \bu_1- \nabla \tbu_1|^2\, \md x \md t + \int_0^{\tau}\int_{\Omega}\,\Big( |\triangle \bd_1 - \triangle \tbd_1|^2+ |\nabla \bd_1 - \nabla \tbd_1|^2\Big)\, \md x \md t\\\label{entr-3}
\leq \,&\, \mathcal{E}_1(0) + \int_0^{\tau}\, \mathcal{R}_1(\rho_1,\bu_1,\bd_1,\tr_1,\tbu_1,\tbd_1)\,\md t,
\end{align}
where
\stepcounter{proof}
\begin{equation}\label{R1}
    \dis \mathcal{R}_1 = \mathcal{R}_1(\rho_1,\bu_1,\bd_1,\tr_1,\tbu_1,\tbd_1):= \mathcal{R}_{1d} + \mathcal{R}_{1c},
\end{equation}
\stepcounter{proof}
\begin{align}\label{R1d}
    \dis \mathcal{R}_{1d}:=\,&\, \int_{\Omega}\, \rho_1 (\tbu_1-\bu_1)\cdot \big((\nabla \tbu_1) (\bu_1-\tbu_1)\big) \, \md x \\\nonumber \,&\,-\int_{\Omega}\,\text{div}\tbu_1\Big(P(\rho_1)-P'(\tr_1)(\rho_1-\tr_1)-P(\tr_1)\Big)\, \md x \\\nonumber
    \,&\,+ \int_{\Omega}\, \frac{(\rho_1-\tr_1)}{\tr_1}\left[\triangle \tbu_1- \text{div}\Big(\nabla \tbd_1 \odot \nabla \tbd_1 - \frac{1}{2}|\nabla \tbd_1|^2 \mathbb{I}_3\Big)\right]\cdot (\tbu_1-\bu_1)\,\md x
\end{align}
and
\stepcounter{proof}
\begin{align}\label{R1c}
    \dis \mathcal{R}_{1c}:=&\, \int_{\Omega}\, \big(\triangle \bd_1 \cdot (\nabla \bd_1)-\triangle \tbd_1 \cdot (\nabla \tbd_1)\big)(\tbu_1-\bu_1)\, \md x  \\\nonumber&\,+ \int_{\Omega}\, (\triangle \bd_1-\triangle \tbd_1) \cdot \big((\nabla \bd_1)\bu_1 - (\nabla \tbd_1)\tbu_1\big)\, \md x  \\\nonumber
&\,-\int_{\Omega}\, (\triangle \bd_1-\triangle \tbd_1) \cdot \big(|\nabla \bd_1|^2\bd_1-|\nabla \tbd_1|^2\tbd_1\big)\, \md x \\\nonumber
&\,+\int_{\Omega}\,(\bd_1 - \tbd_1) \cdot \big(|\nabla \bd_1|^2\bd_1-|\nabla \tbd_1|^2\tbd_1\big)\, \md x \\\nonumber
&\,-\int_{\Omega}\, (\bd_1 - \tbd_1) \cdot \big((\nabla \bd_1)\bu_1 - (\nabla \tbd_1)\tbu_1\big)\, \md x.
\end{align}
Similar to the proof of Theorem \ref{T2.1}, we shall finish the proof of Theorem \ref{T2.2} by applying the relative entropy inequality (\ref{entr-3}) to $\{\tr_1,\tbu_1,\tbd_1\}$, where $\{\tr_1,\tbu_1,\tbd_1\}$ is a classical (smooth) solution of the initial-boundary value problem (\ref{lcd-2}), (\ref{initial}), and (\ref{bou-2}), such that
\begin{equation*}
    \dis \tr_1(0,\cdot)=\rho_0, \quad \tbu_1(0,\cdot)=\bu_0, \quad \tbd_1(0,\cdot)=\bd_0.
\end{equation*}
As a result, the integrals depending on the initial data on the right-hand side of (\ref{entr-3}) vanish,
and we apply a Gronwall type argument to deduce the desired result, namely,
\begin{equation*}
    \dis \rho_1 \equiv \tr_1,\quad \bu_1 \equiv \tbu_1,\quad \bd_1 \equiv \tbd_1.
\end{equation*}
To complete the proof, we have to examine all terms in the remainder (\ref{R1}) and to show that they can be "absorbed" by the left-hand side of (\ref{entr-3}).

From (\ref{initial}), it is clear to see that
\stepcounter{proof}
\begin{align}\nonumber
    \dis \, & \, \left|\int_{\Omega}\, \rho_1 (\tbu_1-\bu_1)\cdot \big((\nabla \tbu_1) (\bu_1-\tbu_1)\big) \, \md x - \int_{\Omega}\,\text{div}\tbu_1\Big(P(\rho_1)-P'(\tr_1)(\rho_1-\tr_1)-P(\tr_1)\Big)\, \md x\right| \\\label{d4}
    \leq \,&\, C \|\nabla \tbu_1 \|_{L^{\infty}(\Omega)}\, \mathcal{E}_1.%\Big([\rho_1,\bu_1,\bd_1]|[\tr_1,\tbu_1,\tbd_1]\Big).
\end{align}
Let
$$\dis \tilde{\mathbf{g}}_1=\tilde{\mathbf{g}}_1(\triangle \tbu_1,\nabla \tbd_1,\triangle \tbd_1)=\triangle \tbu_1- \text{div}\Big(\nabla \tbd_1 \odot \nabla \tbd_1 - \frac{1}{2}|\nabla \tbd_1|^2 \mathbb{I}_3\Big).$$
Obviously, we have
\stepcounter{proof}
\begin{align}\label{d5}
   \dis \,&\, \int_{\Omega}\, \frac{1}{\tr_1}(\rho_1-\tr_1)\, \tilde{\mathbf{g}}_1 \cdot (\tbu_1-\bu_1)\,\md x \\\nonumber
    = \,&\,\int_{\{0 \leq \rho_1 < 2\tr_1\}}\, \frac{1}{\tr_1}(\rho_1-\tr_1)\,\tilde{\mathbf{g}}_1 \cdot (\tbu_1-\bu_1)\,\md x +\int_{\{ \rho_1 \geq 2\tr_1\}}\, \frac{1}{\tr_1}(\rho_1-\tr_1)\, \tilde{\mathbf{g}}_1 \cdot (\tbu_1-\bu_1)\,\md x.
\end{align}
Similar to (\ref{a4}), it follows that
\stepcounter{proof}
\begin{align}\label{d6}
    \dis \left| \int_{\{0 \leq \rho_1 < 2\tr_1\}}\, \frac{1}{\tr_1}(\rho_1-\tr_1)\, \tilde{\mathbf{g}}_1 \cdot (\tbu_1-\bu_1)\,\md x\right| %\\\nonumber
    \leq  C(\delta)\|\frac{\tilde{\mathbf{g}}_1}{\tr_1}\|^2_{L^3(\Omega)}\,\mathcal{E}_1%\Big([\rho_1,\bu_1,\bd_1]|[\tr_1,\tbu_1,\tbd_1]\Big) + \delta \, \|\nabla \bu_1 -\nabla \tbu_1\|_{L^2(\Omega)}^2
\end{align}
for any $\delta > 0$. We proceed in a same way as the derivation of (\ref{a5}) to find that
\stepcounter{proof}
\begin{align}\label{d7}
    \dis  \left|\int_{\{ \rho_1 \geq 2\tr_1\}}\, \frac{1}{\tr_1}(\rho_1-\tr_1)\,\tilde{\mathbf{g}}_1 \cdot (\tbu_1-\bu_1)\,\md x \right|
\dis    \leq  C \, \| \tilde{\mathbf{g}}_1 \|_{L^{\infty}(\Omega)}\,\mathcal{E}_1.
\end{align}

Now we are going to estimate $\mathcal{R}_{1c}$. From (\ref{R1c}), a direct calculation gives
\stepcounter{proof}
\begin{equation}\label{d8}
    \dis \mathcal{R}_{1c} = \mathcal{R}_{1c}^a + \mathcal{R}_{1c}^b %+ \mathcal{R}_{1c}^c
\end{equation}
where
\stepcounter{proof}
\begin{align}\label{d9}
    \dis \mathcal{R}_{1c}^a = \,&\, \int_{\Omega}\, (\triangle \bd_1 - \triangle \tbd_1) \cdot (\nabla \bd_1 - \nabla \tbd_1) \tbu_1 \, \md x \\\nonumber
    \,&\, + \int_{\Omega}\, \triangle \tbd_1 \cdot (\nabla \bd_1 - \nabla \tbd_1) (\tbu_1 - \bu_1) \, \md x \\\nonumber
    \,&\, - \int_{\Omega}\,(|\nabla \tbd_1|+|\nabla \bd_1|) \bd_1 \cdot (\triangle \bd_1 - \triangle \tbd_1)  (|\nabla \bd_1| - |\nabla \tbd_1|)  \, \md x\\\nonumber
    \,&\,+ \int_{\Omega}\, ( \bd_1 -  \tbd_1) \cdot (\nabla \bd_1)(\bu_1 - \tbu_1)  \, \md x
\end{align}
and
\stepcounter{proof}
\begin{align}\label{d10}
    \dis \mathcal{R}_{1c}^b = \,&\, \int_{\Omega}\, ( \bd_1 -  \tbd_1) \cdot (\nabla \bd_1 - \nabla \tbd_1) \tbu_1 \, \md x\\\nonumber
    \,&\, + \int_{\Omega}\,(|\nabla \tbd_1|+|\nabla \bd_1|) \bd_1 \cdot ( \bd_1 -  \tbd_1)  (|\nabla \bd_1| - |\nabla \tbd_1|)  \, \md x \\\nonumber
    \,&\, + \int_{\Omega}\, | \bd_1 -  \tbd_1|^2 |\nabla \tbd_1|^2 \, \md x.
\end{align}
%and
%\stepcounter{proof}
%\begin{align}\label{d11}
%    \dis \mathcal{R}_{1c}^c = \,&\, \int_{\Omega}\, ( \bd_1 -  \tbd_1) \cdot (\nabla \bd_1)(\bu_1 - \tbu_1)  %\, \md x\\\nonumber
%    \,&\, - \int_{\Omega}\,|\nabla \bd_1| \bd_1 \cdot ( \triangle \bd_1 - \triangle \tbd_1)  (|\nabla \bd_1| - %|\nabla \tbd_1|)  \, \md x\\\nonumber
%    \,&\, + \int_{\Omega}\,|\nabla \bd_1| \bd_1 \cdot ( \bd_1 - \tbd_1)  (|\nabla \bd_1| - |\nabla \tbd_1|)  %\, \md x.
%\end{align}
For the first term on the right hand side of (\ref{d9}), we have
\stepcounter{proof}
\begin{align}\label{d12}
    \dis \,&\, \left| \int_{\Omega}\, (\triangle \bd_1 - \triangle \tbd_1) \cdot (\nabla \bd_1 - \nabla \tbd_1) \tbu_1 \, \md x \right| \\\nonumber
    \leq \,&\, \|\tbu_1\|_{L^{\infty}(\Omega)} \|\triangle \bd_1 - \triangle \tbd_1\|_{L^2(\Omega)}\|\nabla \bd_1 - \nabla \tbd_1\|_{L^2(\Omega)}\\\nonumber
    \leq \,&\, \delta\, \|\triangle \bd_1 - \triangle \tbd_1\|^2_{L^2(\Omega)}+ C(\delta) \|\tbu_1\|^2_{L^{\infty}(\Omega)} \|\nabla \bd_1 - \nabla \tbd_1\|^2_{L^2(\Omega)}\\\nonumber
    \leq \,&\, \delta\, \|\triangle \bd_1 - \triangle \tbd_1\|^2_{L^2(\Omega)}+ C(\delta) \|\tbu_1\|^2_{L^{\infty}(\Omega)} \,\mathcal{E}_1,%\Big([\rho_1,\bu_1,\bd_1]|[\tr_1,\tbu_1,\tbd_1]\Big),
\end{align}
for any $\delta > 0$. In a similar way, we also obtain
\stepcounter{proof}
\begin{align}\label{d13}
    \dis \,&\, \left| \int_{\Omega}\, \triangle \tbd_1 \cdot (\nabla \bd_1 - \nabla \tbd_1) (\tbu_1 - \bu_1) \, \md x  \right| \\\nonumber
    \leq \,&\, \delta\, \|\nabla \bu_1 - \nabla \tbu_1\|^2_{L^2(\Omega)}+ C(\delta) \|\triangle \tbd_1\|^2_{L^{\infty}(\Omega)} \,\mathcal{E}_1,%\Big([\rho_1,\bu_1,\bd_1]|[\tr_1,\tbu_1,\tbd_1]\Big),
\end{align}
\stepcounter{proof}
\begin{align}\label{d14}
    \dis \,&\, \left| - \int_{\Omega}\,(|\nabla \tbd_1|+|\nabla \bd_1|) \bd_1 \cdot (\triangle \bd_1 - \triangle \tbd_1)  (|\nabla \bd_1| - |\nabla \tbd_1|)  \, \md x \right| \\\nonumber
    \leq \,&\, \delta\, \|\triangle \bd_1 - \triangle \tbd_1\|^2_{L^2(\Omega)}+ C(\delta)\, \big( \|\nabla \tbd_1\|^2_{L^{\infty}(\Omega)}+\|\nabla \bd_1\|^2_{L^{\infty}(\Omega)}\big)\|\bd_1\|^2_{L^{\infty}(\Omega)} \,\mathcal{E}_1,%\Big([\rho_1,\bu_1,\bd_1]|[\tr_1,\tbu_1,\tbd_1]\Big),
\end{align}
and
\stepcounter{proof}
\begin{align}\label{d18}
    \dis \,&\, \left| \int_{\Omega}\, ( \bd_1 -  \tbd_1) \cdot (\nabla \bd_1)(\bu_1 - \tbu_1)  \, \md x \right| \\\nonumber
    \leq \,&\, \delta \, \| \nabla \bu_1 - \nabla \tbu_1\|^2_{L^{2}(\Omega)} + C(\delta) \,\|\nabla \bd_1\|^2_{L^{\infty}(\Omega)}\,\mathcal{E}_1,%\Big([\rho_1,\bu_1,\bd_1]|[\tr_1,\tbu_1,\tbd_1]\Big).
\end{align}
for any $\delta > 0$.

Finally, Cauchy inequality and the definition of the relative entropy $\mathcal{E}_1$ give
\stepcounter{proof}
\begin{align}\label{d15}
    \dis \,&\, \left| \int_{\Omega}\, ( \bd_1 - \tbd_1) \cdot (\nabla \bd_1 - \nabla \tbd_1) \tbu_1 \, \md x \right| \\\nonumber
    \leq \,&\, \|\tbu_1\|_{L^{\infty}(\Omega)} \| \bd_1 -  \tbd_1\|_{L^2(\Omega)}\|\nabla \bd_1 - \nabla \tbd_1\|_{L^2(\Omega)}\\\nonumber
    \leq \,&\, \|\tbu_1\|_{L^{\infty}(\Omega)}\,\mathcal{E}^{\frac{1}{2}}_1
    \,\mathcal{E}^{\frac{1}{2}}_1\\\nonumber
    \leq \,&\, \|\tbu_1\|_{L^{\infty}(\Omega)}\,\mathcal{E}_1.
\end{align}
Likewise, we get
\stepcounter{proof}
\begin{align}\label{d16}
    \dis \,&\, \left|  \int_{\Omega}\,(|\nabla \tbd_1| + |\nabla \bd_1|) \bd_1 \cdot ( \bd_1 -  \tbd_1)  (|\nabla \bd_1| - |\nabla \tbd_1|)  \, \md x \right| \\\nonumber
    \leq \,&\, \|\bd_1\|_{L^{\infty}(\Omega)} \big(\|\nabla \tbd_1\|_{L^{\infty}(\Omega)}+\|\nabla \bd_1\|_{L^{\infty}(\Omega)}\big) \| \bd_1 -  \tbd_1\|_{L^2(\Omega)}\|\nabla \bd_1 - \nabla \tbd_1\|_{L^2(\Omega)}\\\nonumber
    \leq \,&\, \|\bd_1\|_{L^{\infty}(\Omega)} \big(\|\nabla \tbd_1\|_{L^{\infty}(\Omega)}+\|\nabla \bd_1\|_{L^{\infty}(\Omega)}\big)\,\mathcal{E}_1
\end{align}
and
\stepcounter{proof}
\begin{align}\label{d17}
    \dis   \int_{\Omega}\, | \bd_1 -  \tbd_1|^2 |\nabla \tbd_1|^2 \, \md x
    \leq  \|\nabla \tbd_1\|^2_{L^{\infty}(\Omega)} \| \bd_1 -  \tbd_1\|^2_{L^2(\Omega)}
    \leq   \|\nabla \tbd_1\|^2_{L^{\infty}(\Omega)}\,\mathcal{E}_1.
\end{align}

Moreover, by applying the quasilinear equations of parabolic type estimates (see \cite{Lady68} Chapter VI, Section 3) to the third equation of (\ref{lcd-2}) and a density argument, we have
\stepcounter{proof}
\begin{equation}\label{d19}
    \|\nabla \bd_1\|_{L^{\infty}(\Omega)} \leq C.
\end{equation}

Summing up relations (\ref{R1})--(\ref{d19}) with the definition of finite energy weak solution $\{\rho_1,\bu_1,\bd_1\}$, we conclude from the relative entropy inequality (\ref{entr-3}) that
\[
\dis \mathcal{E}_1\Big([\rho_1,\bu_1,\bd_1]|[\tr_1,\tbu_1,\tbd_1]\Big)(\tau) \leq \int_0^{\tau} \, h_1(t)
\mathcal{E}_1\Big([\rho_1,\bu_1,\bd_1]|[\tr_1,\tbu_1,\tbd_1]\Big)(t)\, \md t,
\]
where $ h_1 \in L^1(0,T)$. Thus, Theorem \ref{T2.1} immediately follows from Gronwall's inequality.

\begin{remark}
\label{R4.1}
For one-dimensional case, Ding {\it et al.} in \cite{DWW11,DLWW12} established the existence of global classical solutions and the existence of global weak solutions. Based on the existence results, one can directly apply Theorem \ref{T2.2} to obtain the weak-strong uniqueness principle for one-dimensional compressible liquid crystal system.
\end{remark}

% Appendix %%%%%%%%%%%%%%%%%%%%%%%%%%%%%%%%%%%%%%%%%%%%%%%%%%%%%%%%%%

%\vspace{1cm}

%\appendix

\vspace{5mm}

\noindent {\large \bf Acknowledgements.}

The authors would like to thank Prof. Song Jiang for helpful discussions and valuable suggestions. This work was supported by NSFC (Grant No. 11171035).
\vspace{2mm}

% ---------------------------------------------------------

\end{document}